\documentstyle[a4]{article}

\begin{document}

\title{A bilinear Airy-estimate with application to gKdV-3}

\author{Axel Gr\"unrock \\Fachbereich Mathematik\\ Bergische Universit\"at - Gesamthochschule Wuppertal \\ Gau{\ss}stra{\ss}e 20 \\ D-42097 Wuppertal \\ Germany \\ e-mail Axel.Gruenrock@math.uni-wuppertal.de}

\date{}

\maketitle

\newcommand{\n}[2]{\mbox{$ \| #1 \| _{ #2} $}}
\newcommand{\XX}[2]{\mbox{$ X_{#1,#2} $}}
\newcommand{\x}{\mbox{$ X_{s,b} $}}
\newcommand{\F}{\mbox{${\cal F}$}}
\newcommand{\q}[2]{\mbox{$ {\| #1 \|}^2_{#2} $}}

\pagestyle{plain}
\rule{\textwidth}{0.5pt}

\newtheorem{lemma}{Lemma}
\newtheorem{kor}{Corollary}
\newtheorem{satz}{Theorem}

\begin{abstract}
The Fourier restriction norm method is used to show local wellposedness for the Cauchy-Problem
\[u_t + u_{xxx} + (u^4)_x=0,\hspace{1cm}u(0)=u_0 \in H^s_x({\bf R}), \,\,\,s>-\frac{1}{6}\]
for the generalized Korteweg-deVries equation of order three, for short gKdV-3. For real valued data $u_0 \in L^2_x({\bf R})$ global wellposedness follows by the conservation of the $L^2$-norm. The main new tool is a bilinear estimate for solutions of the Airy-equation.
\end{abstract}

\vspace{0,5cm}

The purpose of this note is to establish local wellposedness of the Cauchy-Problem
\[u_t + u_{xxx} + (u^4)_x=0,\hspace{1cm}u(0)=u_0 \in H^s_x({\bf R}), \,\,\,s>-\frac{1}{6}\]
for the generalized Korteweg-deVries equation of order three, for short gKdV-3. So far, local wellposedness of this problem is known for data $u_0 \in H^s_x({\bf R})$, $s \ge \frac{1}{12}$. This was shown by Kenig, Ponce and Vega in 1993, see Theorem 2.6 in \cite{KPV93}. Here we extend this result to data  $u_0 \in H^s_x({\bf R})$, $s >- \frac{1}{6}$. A standard scaling argument suggests that this is optimal (up to the endpoint). For real valued data $u_0 \in L^2_x({\bf R})$ we obtain global wellposedness by the conservation of the $L^2$-norm.

By the Fourier restriction norm method introduced in \cite{B93} and further developed in \cite{KPV96} and \cite{GTV97} matters reduce to the proof of the estimate
\begin{equation}\label{1}
\n{\partial_x \prod_{i=1}^4 u_i}{\XX{s}{b'}} \le c \prod_{i=1}^4 \n{u_i}{\XX{s}{b}}
\end{equation}
for suitable values of $s$, $b$ and $b'$. Here the space $\x$ is the completion of the Schwartz class ${\cal S} ({\bf R}^{2})$ with respect to the norm
\[\n{u}{\x} = \left(\int d \xi d\tau <\tau - \xi^3>^{2b}<\xi>^{2s}|\F u(\xi,\tau)|^2 \right) ^{\frac{1}{2}},\]
where $\F$ denotes the Fourier transform in both variables. The main new tool for the proof of (\ref{1}) is the following bilinear Airy-estimate:

\begin{lemma}\label{l1} 
Let $I^s$ denote the Riesz potential of order $-s$ and let $I^s_-(f,g)$ be defined by its Fourier transform (in the space variable):
\[\F_x I_-^s (f,g) (\xi) := \int_{\xi_1+\xi_2=\xi}d\xi_1|\xi_1-\xi_2|^s \F_xf(\xi_1)\F_xg(\xi_2).\]
Then we have
\[\n{I^{\frac{1}{2}}I_-^{\frac{1}{2}}(e^{-t\partial ^3}u_1,e^{-t\partial ^3}u_2)}{L^2_{xt}} \le c \n{u_1}{L^2_{x}}\n{u_2}{L^2_{x}}.\]
\end{lemma}

Proof: We will write for short $\hat{u}$ instead of $\F_x u$ and $\int_* d\xi_1$ for $\int_{\xi_1+\xi_2=\xi}d\xi_1$. Then, using Fourier-Plancherel in the space variable we obtain:

\begin{eqnarray*}
&&  \q{I^{\frac{1}{2}}I_-^{\frac{1}{2}}(e^{-t\partial ^3}u_1,e^{-t\partial ^3}u_2)}{L^2_{xt}} \\
&=& c\int d\xi |\xi | dt \left| \int_* d\xi_1 |\xi_1-\xi_2|^{\frac{1}{2}} e^{it(\xi_1^3 + \xi_2^3)}\hat{u}_1(\xi_1)\hat{u}_2(\xi_2) \right|^2 \\
&=& c\int d\xi |\xi |dt\int_* d\xi_1 d \eta_1 e^{it(\xi_1^3 + \xi_2^3- \eta_1^3 - \eta_2^3)}(|\xi_1-\xi_2||\eta_1-\eta_2|)^{\frac{1}{2}}\prod_{i=1}^2 \hat{u_i}(\xi_i) \overline{\hat{u_i}(\eta_i)} \\
&=& c\int d\xi |\xi |\int_* d\xi_1 d \eta_1 \delta (\eta_1^3 + \eta_2^3 - \xi_1^3 - \xi_2^3)(|\xi_1-\xi_2||\eta_1-\eta_2|)^{\frac{1}{2}}\prod_{i=1}^2 \hat{u_i}(\xi_i) \overline{\hat{u_i}(\eta_i)} \\
&=& c\int d\xi |\xi |\int_* d\xi_1 d \eta_1 \delta (3\xi (\eta_1^2 - \xi_1^2 + \xi(\xi_1-\eta_1)))(|\xi_1-\xi_2||\eta_1-\eta_2|)^{\frac{1}{2}}\prod_{i=1}^2 \hat{u_i}(\xi_i) \overline{\hat{u_i}(\eta_i)}.
\end{eqnarray*}
Now we use $\delta(g(x)) = \sum_n \frac{1}{|g'(x_n)|} \delta(x-x_n)$, where the sum is taken over all simple zeros of $g$, in our case:
\[g(x)= 3 \xi (x^2 + \xi (\xi_1-x)-\xi_1^2)\]
with the zeros $x_1 =\xi_1$ and $x_2 =\xi-\xi_1$, hence $g'(x_1)=3 \xi(2\xi_1-\xi)$ respectively $g'(x_2)=3 \xi (\xi-2\xi_1)$.
So the last expression is equal to
\begin{eqnarray*}
&&c\int d\xi |\xi | \int_* d\xi_1 d \eta_1\frac{1}{|\xi||2\xi_1-\xi|}\delta(\eta_1-\xi_1)(|\xi_1-\xi_2||\eta_1-\eta_2|)^{\frac{1}{2}}\prod_{i=1}^2 \hat{u_i}(\xi_i) \overline{\hat{u_i}(\eta_i)}\\
&+& c\int d\xi |\xi | \int_* d\xi_1 d \eta_1\frac{1}{|\xi||2\xi_1-\xi|}\delta(\eta_1-(\xi-\xi_1))(|\xi_1-\xi_2||\eta_1-\eta_2|)^{\frac{1}{2}}\prod_{i=1}^2 \hat{u_i}(\xi_i) \overline{\hat{u_i}(\eta_i)}\\
&=&c\int d\xi \int_* d\xi_1 \prod_{i=1}^2|\hat{u_i}(\xi_i)|^2 + c\int d\xi \int_* d\xi_1\hat{u}_1(\xi_1)\overline{\hat{u}_1}(\xi_2)\hat{u}_2(\xi_2)\overline{\hat{u}_2}(\xi_1)\\
&\le& c (\prod_{i=1}^2 \q{u_i}{L^2_x} + \q{\hat{u}_1\hat{u}_2}{L^1_{\xi}}) \le c \prod_{i=1}^2 \q{u_i}{L^2_x}.
\end{eqnarray*}
$\hfill \Box$

Arguing as in the proof of Lemma 2.3 in \cite{GTV97} we get the following

\begin{kor}\label{k1} Let $b > \frac{1}{2}$. Then the following estimate holds true:
\[\n{I^{\frac{1}{2}}I_-^{\frac{1}{2}} (u,v)}{L^2_{xt}}\le c \n{u}{\XX{0}{b}}\n{v}{\XX{0}{b}}\]
\end{kor}

\vspace{0.5cm}

In the next Lemma, some well known Strichartz type estimates for the Airy equation are gathered in terms of $\x$-norms:

\begin{lemma}\label{l2} For $b > \frac{1}{2}$ the following estimates are valid:
\begin{itemize}
\item[i)] $\n{u}{L_t^p(H_x^{s,q})} \le c \n{u}{\XX{0}{b}}$, whenever $0\le s=  \frac{1}{p} \le \frac{1}{4}$ and $\frac{1}{q}=\frac{1}{2}-\frac{2}{p}$,
\item[ii)] $\n{u}{L_t^p(L_x^{q})} \le c \n{u}{\XX{0}{b}}$, whenever $0<\frac{1}{q}=\frac{1}{2}-\frac{3}{p} \le \frac{1}{2}$.
\end{itemize}
\end{lemma}

Quotation/proof: Theorem 2.1 in \cite{KPV91} gives in the case of the Airy-equation
\[\n{e^{-t\partial ^3} u_0}{L_t^p(\dot{H}_x^{s,q})} \le c \n{u_0}{L_x^2},\]
provided $0\le s=  \frac{1}{p} \le \frac{1}{4}$ and $\frac{1}{q}=\frac{1}{2}-\frac{2}{p}$. Now Lemma 2.3 in \cite{GTV97} is applied to obtain
\begin{equation}\label{2}
\n{u}{L_t^p(\dot{H}_x^{s,q})} \le c \n{u}{\XX{0}{b}}, \hspace{1cm}b>\frac{1}{2}
\end{equation}
for the same values of $s$, $p$ and $q$. From this ii) follows by Sobolev's embedding theorem (in the space variable). Especially we have
\[\n{u}{L^8_{xt}} \le c \n{u}{\XX{0}{b}}, \hspace{1cm}b>\frac{1}{2}
,\]
which, interpolated with the trivial case, gives
\[\n{u}{L^4_{xt}} \le c \n{u}{\XX{0}{b}}, \hspace{1cm}b>\frac{1}{3}.\]
Now let us see how to replace $\dot{H}_x^{s,q}$ by $H_x^{s,q}$ in (\ref{2}) in the endpoint case, i. e. $s=\frac{1}{p}=\frac{1}{4}$, $q=\infty$: Using the projections $p=\F _x^{-1}\chi_{\{|\xi|\le 1\}}\F _x$ and $P=Id - p$ we have
\[\n{u}{L_t^4(H_x^{\frac{1}{4},\infty})} \le \n{Pu}{L_t^4(H_x^{\frac{1}{4},\infty})}+\n{pu}{L_t^4(H_x^{\frac{1}{4},\infty})}=:I+II.\]
For $I$ we use (\ref{2}) to obtain
\[I \le c\n{I^{-\frac{1}{4}}J^{\frac{1}{4}}Pu}{\XX{0}{b}} \le c \n{u}{\XX{0}{b}},\]
while for $II$ by Sobolev's embedding theorem we get
\[II \le c \n{pu}{L_t^4(H_x^{\frac{1}{2}+,4})} \le c \n{pu}{\XX{\frac{1}{2}+}{b}}\le c \n{u}{\XX{0}{b}}.\]
This gives i) in the endpoint case, from which the general case follows by interpolation with Sobolev's embedding theorem (in the time variable).
$\hfill \Box$

$Remark$: The endpoint case in ii) is also valid - see e. g. Lemma 3.29 in \cite{KPV93} - but we shall not make use of this here.

Now we are prepared to prove the crucial nonlinear estimate:

\begin{satz} For $0\ge s > -\frac{1}{6}$, $-\frac{1}{2}<b'<s-\frac{1}{3}$ and $b>\frac{1}{2}$ the estimate (\ref{1}) is valid.
\end{satz}

Proof: Writing $f_i(\xi, \tau)= <\tau - \xi^3>^{b}<\xi>^s \F u_i(\xi, \tau)$, $1 \le i \le 4$, we have
\[\n{ \partial_x \prod_{i=1}^4 u_i}{\XX{s}{b'}}=c \n{<\!\!\tau \!\!-\!\! \xi^3\!\!>^{b'} <\!\!\xi\!\!>^{s}|\xi|\int d \nu  \prod_{i=1}^4 <\!\!\tau_i\!\! -\!\!\xi_i^3\!\!>^{-b} <\!\!\xi_i\!\!>^{-s}f_i(\xi_i,\tau_i)}{L^2_{\xi,\tau}},\]
where $d \nu = d\xi_1..d\xi_{3} d\tau_1.. d \tau_{3}$ and $\sum_{i=1}^4 (\xi_i,\tau_i) = (\xi, \tau)$. Now the domain of integration is devided into the regions $A$, $B$ and $C=(A \cup B)^c$, where in $A$ we assume $|\xi_{max}| \le c$. (Here $\xi_{max}$ is defined by $|\xi_{max}| = \max_{i=1}^4 |\xi_i|$, similarly $\xi_{min}$.) Then for the region $A$ we have the upper bound
\begin{eqnarray*}
&& c\n{\int d \nu  \prod_{i=1}^4 <\!\!\tau_i\!\! -\!\!\xi_i^3\!\!>^{-b}f_i(\xi_i,\tau_i)}{L^2_{\xi,\tau}} \\ &=& c \n{\prod_{i=1}^4 J^s u_i}{L^2_{x,t}} 
\le c \prod_{i=1}^4 \n{J^s u_i}{L^8_{x,t}}  \le  c \prod_{i=1}^4 \n{u_i}{\x},
\end{eqnarray*}
where in the last step Lemma \ref{l2} , part ii), with $p=q=8$ was applied.

Besides $|\xi_{max}| \ge c$ ($\Rightarrow <\xi_{max}> \le c |\xi_{max}|$) we shall assume for the region $B$ that
\begin{itemize}
\item[i)] $|\xi_{min}| \le 0.99 |\xi_{max}|$ or
\item[ii)] $|\xi_{min}| > 0.99 |\xi_{max}|$, and there are exactly two indices $i \in \{1,2,3,4 \}$ with $\xi_i > 0$.
\end{itemize}
Then the region $B$ can be splitted again into a finite number of subregions, so that for any of these subregions there exists a permutation $\pi$ of $\{1,2,3,4 \}$ with
\[|\xi|<\xi>^s \prod_{i=1}^4<\xi_i>^{-s} \le c |\xi_{\pi (1)}+\xi_{\pi (2)}|^{\frac{1}{2}}|\xi_{\pi (1)}-\xi_{\pi (2)}|^{\frac{1}{2}}<\xi_{\pi (3)}>^{-\frac{3s}{2}}<\xi_{\pi (4)}>^{-\frac{3s}{2}}.\]
Assume $\pi = id$ for the sake of simplicity now. Then we get the upper bound
\begin{eqnarray*}
 \n{\!<\!\!\tau \!\!-\!\! \xi^3\!\!>^{b'}\!\int\! d \nu|\xi_{1}+\xi_{2}|^{\frac{1}{2}}|\xi_{1}-\xi_{2}|^{\frac{1}{2}}\!<\!\!\xi_{3}\!\!>\!\!^{-\frac{3s}{2}}\!<\!\!\xi_{4}\!\!>\!\!^{-\frac{3s}{2}}\prod_{i=1}^4\!<\!\!\tau_i\!\! -\!\!\xi_i^3\!\!>^{-b}\!f_i(\xi_i,\tau_i)}{L^2_{\xi,\tau}} \\
= c \n{(I^{\frac{1}{2}}I_-^{\frac{1}{2}}(J^su_1,J^su_2))(J^{-\frac{s}{2}}u_3)(J^{-\frac{s}{2}}u_4)}{\XX{0}{b'}}\hspace{3cm}
\end{eqnarray*}
To estimate the latter expression, we fix some Sobolev- and H\"olderexponents:
\begin{itemize}
\item[i)] $\frac{1}{q_0}=\frac{1}{2}-b'$ so that $L_t^{q_0}(L_x^2) \subset \XX{0}{b'}$,
\item[ii)] $\frac{2}{p}=\frac{1}{q_0}-\frac{1}{2}=-b'$,
\item[iii)] $\frac{1}{q}=\frac{1}{2}-\frac{2}{p}=\frac{1}{2}+b'$ so that by Lemma \ref{l2} $\n{J^{\frac{1}{p}}u}{L_t^{p}(L_x^q)} \le c \n{u}{\XX{0}{b}}$,
\item[iv)] $\epsilon = \frac{1}{p}+\frac{3s}{2} > \frac{1}{q}$ (since $s>\frac{1}{3} + b'$) so that $H_x^{\epsilon,q} \subset L_x^{\infty}$.
\end{itemize}
Now we have
\begin{eqnarray*}
&&\n{(I^{\frac{1}{2}}I_-^{\frac{1}{2}}(J^su_1,J^su_2))(J^{-\frac{s}{2}}u_3)(J^{-\frac{s}{2}}u_4)}{\XX{0}{b'}}\\
& \le & c \n{(I^{\frac{1}{2}}I_-^{\frac{1}{2}}(J^su_1,J^su_2))(J^{-\frac{s}{2}}u_3)(J^{-\frac{s}{2}}u_4)}{L_t^{q_0}(L_x^2)}\\
& \le & c \n{I^{\frac{1}{2}}I_-^{\frac{1}{2}}(J^su_1,J^su_2)}{L^2_{xt}} \n{J^{-\frac{s}{2}}u_3}{L_t^p(L_x^{\infty})}\n{J^{-\frac{s}{2}}u_4}{L_t^p(L_x^{\infty})}.
\end{eqnarray*}
Now by Lemma \ref{1} the first factor can be controlled by $c\n{u_1}{\x}\n{u_2}{\x}$, while for the second we have the upper bound
\[c \n{J^{-\frac{3s}{2}+\epsilon}J^su_3}{L_t^p(L_x^{q})}= c\n{J^{\frac{1}{p}}J^su_3}{L_t^p(L_x^{q})}\le c \n{u_3}{\XX{s}{b}}.\]
The third factor can be treated in presicely the same way. So for the contributions of the region $B$ we have obtained the desired bound.

Finally we consider the remaining region $C$: Here the $|\xi_i|$, $1\le i \le 4$, are all very close together and $\ge c <\xi_i>$. Moreover, at least three of the variables $\xi_i$ have the same sign. Thus for the quantity $c.q.$ controlled by the expressions $<\tau - \xi ^3>$, $<\tau_i - \xi_i ^3>$, $1\le i \le 4$, we have in this region:
\[c.q.:=|\xi ^3 - \sum_{i=1}^4 \xi_i ^3| \ge c \sum_{i=1}^4 <\xi_i >^3 \ge c <\xi >^3\]
and hence, since $s>\frac{1}{3}+b'$ is assumed,
\[|\xi|<\xi>^s \prod_{i=1}^4<\xi_i>^{-s} \le c(<\tau - \xi ^3>^{-b'}+ \sum_{i=1}^4<\tau_i - \xi_i ^3>^{-b'}\chi_{C_i}),\]
where in the subregion $C_i$, $1\le i \le 4$,  the expression $<\tau_i - \xi_i ^3>$ is dominant. The first contribution can be estimated by
\begin{eqnarray*}
&& c\n{\int d \nu  \prod_{i=1}^4 <\!\!\tau_i\!\! -\!\!\xi_i^3\!\!>^{-b}f_i(\xi_i,\tau_i)}{L^2_{\xi,\tau}} \\ &=& c \n{\prod_{i=1}^4 J^s u_i}{L^2_{x,t}} 
\le c \prod_{i=1}^4 \n{J^s u_i}{L^8_{x,t}}  \le  c \prod_{i=1}^4 \n{u_i}{\x},
\end{eqnarray*}
where we have used Lemma \ref{l2}, part ii). For the contribution of the subregion $C_1$ we take into account that $<\tau_1 - \xi_1 ^3>=\max\{<\tau - \xi ^3>,<\tau_i - \xi_i ^3>,1\le i \le 4\}$, which gives
\[<\tau - \xi ^3>^{b+b'}|\xi|<\xi>^s \prod_{i=1}^4<\xi_i>^{-s} \le c <\tau_1 - \xi_1 ^3>^{b}.\]
So, for this contribution we get the upper bound
\begin{eqnarray*}
&& c\n{<\tau - \xi ^3>^{-b}\int d \nu <\tau_1 - \xi_1 ^3>^{b} \prod_{i=1}^4 <\!\!\tau_i\!\! -\!\!\xi_i^3\!\!>^{-b}f_i(\xi_i,\tau_i)}{L^2_{\xi,\tau}} \\
&\le & c \n{\F ^{-1}f_1 \prod_{i=2}^4 J^s u_i}{\XX{0}{-b}} \le c\n{\F ^{-1}f_1 \prod_{i=2}^4 J^s u_i}{L^{\frac{8}{7}}_{xt}} \\
& \le & c\n{\F ^{-1}f_1}{L^2_{xt}}\prod_{i=2}^4 \n{J^s u_i}{L^8_{x,t}} \le  c \prod_{i=1}^4 \n{u_i}{\x}.
\end{eqnarray*}
Here we have used the dual version of the $L^8$-Strichartz estimate, H\"older and the estimate itself. For the remaining subregions $C_i$ the same argument applies.
$\hfill \Box$

\begin{kor} For $s \ge 0$, $-\frac{1}{2}<b'<-\frac{1}{3}$ and $b>\frac{1}{2}$ the estimate (\ref{1}) holds true.
\end{kor}

Proof: For $s=0$ this is contained in the above theorem, while for $s>0$ one only has to use $<\xi> \le c \prod_{i=1}^4<\xi_i>$.
$\hfill \Box$

\vspace{0,5cm}

{\bf{Acknowledgement:}} I want to thank Professor Hartmut Pecher for numerous helpful conversations.

\end{document}